%% file: hantke_mwy_ana.tex
\documentclass[12pt,twoside,titlepage,openright,headsepline]{article}

\usepackage[english]{babel}
\usepackage{multicol,longtable}
\usepackage{lineno,hyperref,amsmath,amsfonts,graphics}
\usepackage{amsfonts,exscale,amssymb,stmaryrd,amsthm}
\usepackage{color}
\usepackage{comment}
\usepackage[latin1]{inputenc}
\usepackage{hvfloat}
\usepackage[T1]{fontenc}
\usepackage{graphicx,epsfig}
\usepackage{epstopdf}
\usepackage{bibgerm,ngerman}
\usepackage{textcomp}
\usepackage{cite}
\usepackage{pifont} 
\usepackage{todonotes}

\DeclareMathOperator{\dv}{div}
\DeclareMathOperator{\D}{d}

\newtheorem{teo}{Theorem}
\newtheorem{rem}{Remark}

\begin{document}


\title {The Riemann problem for a two-phase mixture hyperbolic system with phase function and multi-component equation of state}

\author{
Maren Hantke\thanks{ \baselineskip=.5\baselineskip Institut f\"ur Mathematik, Martin-Luther-Universit\"at Halle - Wittenberg, 
D-06099 Halle (Saale), Germany. Email: maren.hantke@mathematik.uni-halle.de},
Gerald Warnecke\thanks{ \baselineskip=.5\baselineskip Institute of Analysis and  Numerics,
Otto-von-Guericke-University Magdeburg, PSF 4120, D--39016 Magdeburg,
Germany.  Email: warnecke@ovgu.de},
Christoph Matern\thanks{ \baselineskip=.5\baselineskip Mathematisches Institut der Heinrich-Heine-Universit\"at D\"usseldorf, 
	Universit\"atsstr.\ 1, D--40225 D\"usseldorf, Germany. Email: christoph.matern@hhu.de},
Hazem Yaghi\thanks{ \baselineskip=.5\baselineskip Institute of Analysis and  Numerics,
Otto-von-Guericke-University Magdeburg, PSF 4120, D--39016 Magdeburg,
Germany.  Email: hazem.yaghi@ovgu.de}
}

\date {\today }

\maketitle

\begin{abstract}
In this paper a 
hyperbolic system of partial differential equations for two-phase mixture flows with $N$ components
is studied. It is derived from a more complicated model involving diffusion and exchange terms. Important features of the model
are the assumption of isothermal flow, the use of a phase field function to distinguish the phases 
and a mixture equation of state involving the phase field function as well as an affine relation between partial densities and partial pressures in
the liquid phase. This complicates the analysis. A complete solution of the Riemann
initial value problem is given. Some interesting examples are suggested as bench marks for numerical schemes. 
\end{abstract}

\input sec1.tex 
\input sec2.tex 
\input sec3.tex 
\input sec4.tex 
\input sec5.tex 

\section{Conclusion and outlook}

In this paper we presented the exact Riemann solution for the hyperbolic 2-phase multi-component model \eqref{hyp}. Existence and uniqueness of this solution are proven. We have given two examples that may be used as benchmark tests in numerical simulations. Future work will focus on an appropriate Riemann solver.

\bibliographystyle{abbrv}
\bibliography{bibfiles_hmwy_ana}

\end{document}

%% file: sec1.tex
\section{Introduction}

Two-phase flows of liquid and gas with multiple components have a wide range of applications in nature and chemical engineering.
Examples are atmospheric flow or bubble reactors. It is a challenge to
model the interactions of the fluids, especially the exchange of mass and energy due to phase transitions and chemical 
reactions. Several types of models with different advantages and disadvantages are available in the literature. 

Models of Baer-Nunziato type \cite{baernun} typically require a large number of equations that increase numerical cost substantially. 
Further, these models are usually not in divergence form. Accordingly special attention has to be given 
to their discretization. Moreover, the 
form of the exchange terms is not known. They have to be derived for the special situation at hand.
For more details see H\`erard \cite{herard} and M\"uller et al.\ \cite{simapa}. 

In 2007 Romenski et al.\ \cite{romenski} introduced the similar symmetric hyperbolic and thermodynamically compatible (SHTC) two-phase flow model. 
Although the volume fraction is a variable of the system the model is in divergence form. Therefore the model seems to be very interesting from a 
mathematical point of view even though a conservation law for relative velocity should be discussed extensively. Recently all characteristic fields of the system and all possible wave phenomena were discussed, see Thein et al.\ \cite{theinRD}, but the full Riemann solution is still not available.

Sharp interface models need only a smaller number of equations. 
Interesting analytical results are available in Hantke et al.\ \cite{hdw}. Here the Riemann problem for the isothermal Euler equations with phase transitions was completely discussed. Here mass transfer is modeled by a kinetic relation.
To solve such systems numerically, the interface has to be resolved more or less
exactly. Accordingly, either a very fine grid resolution is required, or the computations have to be performed on a 
moving mesh or one has to track the interfaces on an additional mesh. This can become quite complicated in 
higher space dimensions,
see for instance Chalons et al.\ \cite{chalons1} and \cite{chalons2}.
In \cite{chalons1} a conservative finite volume method was developed to approximate weak 1d-solutions of conservation laws with phase boundaries. 
This method was generalized in \cite{chalons2} for 2d-computations and is able to exactly resolve planar phase interfaces. Further interesting results on this topic can be found in Schleper \cite{schleper} or Fechter et al.\ \cite{fechter}.

To overcome the disadvantages mentioned above one can consider phase field models. 
For these a phase field parameter is introduced as a function of time and space. This parameter takes two distinct 
values to indicate the local phases. It smoothly 
changes at the interfaces which are modeled as small zones of finite width.
There is a growing interest in the use of phase field models for this purpose. The phase field can be compared
to the volume fraction in Baer-Nunziato type models. It does not appear as often in the equations but in the
equation of state. This lead to specific challenges for the analysis below as well as for numerical computations, 
see Hantke et al.\ \cite{hmwy}.

The phase field model derived by Dreyer, Giesselmann and Kraus \cite{dgk} is in the focus of this paper. We consider
the conservation law sub-part of this model by neglecting dissipative terms and source functions. The aim of this paper is to
give a full analysis of the Riemann initial value problem for this hyperbolic system of conservation laws. The resulting
system looks rather harmless. But, the challenge comes from the very complicated nature of the equation of state that is needed,
see \eqref{pressure}. The partial pressures are affine functions in the liquid case and there is a dependence on the phase field function 
$\chi$. A good knowledge of the hyperbolic sub-part of the full model is important for numerical computations since it needs
the most careful treatment in order to achieve accurate solutions.

The model is introduced in Section 2. Then we discuss the main mathematical properties of the system in Section 3.
The various cases that may occur in the solutions to Riemann problems are considered in the fourth section. 
It turns out that for the special case of a single component flow with pure phases the result obtained here coincides with a result presented in Hantke et al.\ \cite{hdw}.
In Section 5 we give
some interesting examples that may serve as benchmarks for numerical schemes.
Finally, we draw some conclusions and give an outlook on future work.

%% file: sec2.tex
\section{The model}

Let us consider a mixture of $N$ components. For its description we introduce
the mixture density $\rho$ as well as the partial densities $\rho_\alpha$ for $\alpha=1,\dots,N$
of the $N$ components. They are related by
\begin{equation}
\label{mixden}
\rho =\sum_{\alpha =1}^N \rho_\alpha.
\end{equation} 
Therefore, only $N$ of the densities have to be determined by partial differential equations. 
Next let us denote by $m_\alpha$ the atomic mass of 
component $\alpha$,  and $\mu_{\alpha}$ the 
chemical potential of constituent $\alpha$ and $\gamma_{\alpha}^i$ the stoichiometric coefficients of $N_R$ possible 
chemical reactions, $A^i$ the affinities and $M_{\alpha\beta},M_r^i, M_p$ the mobilities.

Further quantities relevant for the model are
the mixture velocity $\mathbf{v}$, the pressure $p$, temperature $T$, the free energy density of the mixture $\rho\psi$,
the Boltzmann constant $k$ 
and the phase field $\chi$ in $[0,T)\times\Omega$, $\Omega\subset\mathbb{R}^d$.
The phase field quantity $\chi$ denotes the present phase. It takes values in the interval $[-1,1]$. Here  we choose $\chi=1$ to indicate
pure liquid, while $\chi=-1$ indicates a pure vapor phase.

The compressible $N$-component two-phase model introduced by Dreyer et al.\ \cite{dgk} is given by the following system of 
partial differential equations for $\alpha = 1\ldots ,N-1$
\begin{eqnarray*}
	\partial_t\rho+\dv(\rho\mathbf{v}) &=& 0 
	\\
	\partial_t\rho_{\alpha}+\dv(\rho_{\alpha}\mathbf{v})-
	\dv\left(\sum_{\beta=1}^{N-1}M_{\alpha\beta}\nabla(\mu_{\alpha}-\mu_N)\right)
	=&\displaystyle
	\sum_{i=1}^{N_R}&\gamma_{\alpha}^{i}m_{\alpha}M_r^{i}\left(1-\exp\left(\frac{A^i}{kT}\right)\right)
	\\
	\partial_t(\rho\mathbf{v})+\dv(\rho\mathbf{v}\otimes\mathbf{v})+\nabla p 
	+\dv(\gamma\nabla\chi\otimes\nabla\chi-\mathbf{\sigma}_{NS}) &=& 0
	\\
	\rho\partial_t\chi+\rho\mathbf{v}\cdot\nabla\chi
	&=&
	-M_p\left(\frac{\partial \rho\psi}{\partial\chi}-\gamma\Delta\chi\right)
\end{eqnarray*}

This model considers isothermal chemically reacting viscous liquid-vapor flows of $N$ constituents with phase transitions in 
$d$ dimensions. To solve such systems numerically one typically uses splitting methods. This means, that the system 
is split into two sub-problems. On the one hand one considers the flow part of the system, on the other the reacting part by 
integrating the source terms. In this work we will focus on the flow part. Accordingly in the following we restrict ourselves 
to the 1-$d$ homogeneous subsystem of first order terms that is given by
\begin{subequations}\label{hyp}
	\begin{eqnarray}
		\partial_t\rho_{\alpha}+\partial_x(\rho_{\alpha}v)&=&0 \qquad \alpha=1,\dots,N\\
		\partial_t(\rho v)+\partial_x(\rho v^2+p)&=&0\\
		\label{chie}
		\rho\partial_t\chi+\rho v\partial_x\chi&=&0\,.
	\end{eqnarray}
\end{subequations}
Here we replaced the balance of mass for the mixture density $\rho$ by the mass balance 
equation for constituent $N$. In terms of equations and jump conditions this is a fully equivalent replacement.
We refer to this model as the hyperbolic system \eqref{hyp}.

\paragraph{Equation of state.}

The pressure $p$ is a constitutive quantity that is related to the phase field variable $\chi$ and the partial densities 
$\rho_{\alpha}$, $\alpha=1,\dots,N$ of the components by an equation of state $p=p(\chi,\rho_1,\dots,\rho_N)$.
This equation of state is derived from the free energy density 
\[
\rho\psi=W(\chi)+h(\chi)\rho\psi_{L}(\rho_1,\dots,\rho_N)+(1-h(\chi))\rho\psi_V(\rho_1,\dots,\rho_N),
\] 
of the mixture,
where $W(\chi)=w_0(\chi-1)^2(\chi+1)^2$ denotes the double well potential and $\rho\psi_L$ and $\rho\psi_V$ are the free 
energy density functions of the liquid and the vapor phases, resp.

There is some degree of freedom in choosing $h$ and the free energy densities. We follow Dreyer et al.\ \cite{dgk} and choose 
$h\,:\,\mathbb{R}\to[0,1]$ to be 
\[
h(\chi)=\begin{cases}
	0& \chi\le-1\\
	(-\frac{1}{4}\chi+\frac{1}{2})(\chi+1)^2 & -1< \chi< 1\\
	1 & \chi\ge1
\end{cases}\,.
\]
This is the simplest smooth interpolation function satisfying
\[
h(-1)=0\, \quad h(1)=1 \quad \mbox{and} \quad h'(\chi)=0 \quad \mbox{for} \quad |\chi|\ge1.
\]

For the free energy density functions
\[
\rho\psi_L=\sum_{\alpha=1}^N\rho_{\alpha}\psi_{L\alpha} \quad \mbox{and} \quad 
\rho\psi_V=\sum_{\alpha=1}^N\rho_{\alpha}\psi_{V\alpha}
\]
we proceed as follows. We use the partial free energy density functions $\psi_{k\alpha}$, $k\in\{L,V\}$, $\alpha=1,\dots,N$  
such that in pure phases in the single 
component case we end up with the stiffened gas law. 
In the isothermal case this simply means that $\rho_\alpha\psi_{k\alpha}$ with $k\in\{V,L\}$ are chosen to be
\[
\rho_\alpha\psi_{k\alpha}=a_{k\alpha}^2\rho_\alpha\ln\frac{\rho_{\alpha}}{\rho_{k\alpha,ref}}-d_{k\alpha}+\rho_{\alpha}\psi_{k\alpha,ref}\,.
\]
It is known that the chemical potentials $\mu_{k\alpha}$ and the partial pressures $p_{k\alpha}$ are defined by
\[
\mu_{k\alpha}=\frac{\partial\rho_\alpha\psi_{k\alpha}}{\partial\rho_\alpha} \qquad \mbox{and} \qquad p_{k\alpha}=-\rho_\alpha\psi_{k\alpha}+\rho_{\alpha}\mu_{k\alpha} \qquad \mbox{resp.}
\]
For the partial pressures $p_{k\alpha}$ one can easily verify that $p_{k\alpha}(\rho_\alpha)=a^2_{k\alpha}\rho_{\alpha}+d_{k\alpha}$. Here $a_{k\alpha}$ is the isothermal sound speed of 
component $\alpha$ in phase $k\in\{L,V\}$. 
The parameter $d_{k\alpha}$ equals zero for ideal gases.

The mixture pressure $p$ is defined by
\[
p=-\rho\psi+\sum_{\alpha=1}^N\frac{\partial \rho\psi}{\partial \rho_{\alpha}}.
\] 
This leads to
\begin{equation}\label{pressure}
	p(\chi,\rho_1,\dots,\rho_N)=-W(\chi)+
	h(\chi)\sum_{\alpha=1}^N p_{L\alpha}(\rho_{\alpha})+(1-h(\chi))\sum_{\alpha=1}^N
	p_{V\alpha}(\rho_{\alpha}).
\end{equation}
For more details on this and the general case see Dreyer and Bothe \cite{db} as well as Hantke and M\"uller \cite{hm}.

The double well function moderates the mass transfer between the phases due to condensation and evaporation. The interpolation 
function $h$ allows us to describe mixtures of the vapor and the liquid phase.

\paragraph{Discontinuities, jump conditions.}

The equation for the phase function \eqref{chie} is not conservative. But by summing the equations for
the partial densities as in \eqref{mixden} we obtain the equation for the mixture mass density $\rho_t + (\rho u)_x =0$ from
the original model. 
We can multiply this equation by $\chi$, the equation \eqref{chie} by $\rho$ and add the equations. Then we obtain
the conservative equation $(\rho\chi)_t + (\rho v\chi)_x =0$ that is a smoothly equivalent
replacement for the the last equation of the system \eqref{hyp}. Though there seems to be no physical
interpretation for the conserved quantity $\rho\chi$, it may be used to find jump conditions for the advective equation \eqref{chie}.

For a discontinuity traveling at a speed $s=\frac{\D x}{\D t}$ we use for the states on the left the subscript minus $-$ and on the right
plus $+$.
The conservative equation for $\rho\chi$ has the jump condition 
\begin{equation}
\label{jump_chi}
s[\rho_-\chi_- -\rho_+\chi_+] =[\rho_-v_-\chi_- -\rho_+v_+\chi_+].
\end{equation}
The jump condition for the mixture density $s[\rho_--\rho_+] =[\rho_-v_--\rho_+v_+]$ can be solved for $\rho_-v_-$
and the result inserted into \eqref{jump_chi}. Due to cancellations of terms one obtains 
the jump condition $\rho_+(s-v_+)[\chi_--\chi_+]=0$.
Analogously, doing the same with $\rho_+u_+$ leads to a second jump condition
$\rho_-(s-v_-)[\chi_--\chi_+]=0$. Together we obtain the two jump conditions
\begin{equation}
\label{jump_chi_mod}
\rho_\pm (s-v_\pm )[\chi_--\chi_+]=0.
\end{equation}
We may assume that there is no vacuum state, i.e.\ $\rho>0$ everywhere.
The jump conditions \eqref{jump_chi_mod} then imply for any discontinuity in $\chi$, i.e.\ satisfying $[\chi_--\chi_+]\ne 0$,  
that $s=v_-=v_+$. Therefore, the phase function $\chi$ only has jump discontinuities along contact discontinuities
that travel with the velocity $v$ of the flow.  
The advective equation $\chi_t + v\chi_x =0$ is linear in $\chi$
for given $v$. So the jump condition is exactly what one would expect for such an advection equation.

\paragraph{Special case.}
For the single component case $N=1$ and pure phases, i.e.\ $\chi\in\{-1,1\}$ the resulting system is identical to the isothermal 
Euler. Existence and uniqueness results for two-phase Riemann problems were completely discussed in \cite{hdw}.

%% file: sec3.tex
\section{Mathematical properties of the system}\label{prop}

\paragraph{Riemann problem.}
In the following we want to study the Riemann problem for the hyperbolic system. Such a Riemann problem is given by the 
balance equations \eqref{hyp}, the equation of state \eqref{pressure} and the Riemann initial data
\begin{align}\label{ind}
	\chi(0,x)=\begin{cases}
		\chi_{-} \,,& \! x<0\\
		\chi_{+} \,,& \! x>0
	\end{cases}	
	\quad,\quad
	\rho(0,x)=\begin{cases}
		\rho_{\alpha-}\,, &\! x<0\\
		\rho_{\alpha+}\,, &\! x>0
	\end{cases}	
	\quad,\quad
	v(0,x)=\begin{cases}
		v_{-} \,,& \! x<0\\
		v_{+} \,,& \! x>0
	\end{cases}	
\end{align}
for $\alpha=1,\dots,N$.
These initial data are merged in the vectors $\mathbf{W}_-=(\chi_-,\rho_{1-},\dots,\rho_{N-},v_-)^T$ and 
$\mathbf{W}_+=(\chi_+,\rho_{1+},\dots,\rho_{N+},v_+)^T$.

The solution of the Riemann problem consists of constant states that are separated by classical waves, i.e.\ shock waves, 
rarefactions and contacts.

\paragraph{Jacobian matrix.}
In order to determine the structure of the Riemann solution we rewrite the system \eqref{hyp} in its quasilinear form for the 
primitive variables $\chi, \rho_\alpha,v$ with $\alpha=1,\dots, N$.

We have $p=p(\chi,\rho_1,\dots,\rho_N)$, this gives
\[
\partial_x p= \sum_{\alpha=1}^N \partial_{\rho_\alpha}p \partial_x \rho_\alpha
+\partial_\chi p \partial_x \chi\,.
\]
From the equation of state \eqref{pressure} we determine
\begin{equation}\label{Arel}
	\partial_{\rho_\alpha}p=h(\chi)a_{L\alpha}^2+(1-h(\chi))a_{V\alpha}^2=:A^2_{\alpha}
\end{equation}
and
\begin{equation}\label{Brel}
	\partial_\chi p=-W'(\chi)+h'(\chi)\sum_{\alpha=1}^N(a_{L\alpha}^2\rho_\alpha+d_{L\alpha})-h'(\chi)
	\sum_{\alpha=1}^Na_{V\alpha}^2\rho_\alpha=:B\,.
\end{equation}
Here we should keep in mind, that abbreviations $A_\alpha$ depend on $\chi$, while $B$ depends on $\chi$ and the partial densities 
$\rho_\alpha$.

Using the above notations we obtain the quasilinear form of \eqref{hyp} that reads
\begin{subequations}
	\begin{align}
		\partial_t\chi+v\partial_x\chi &= 0\\
		\partial_t\rho_\alpha+v\partial_x\rho_\alpha+\rho_{\alpha}\partial_x v &=0 \quad \alpha=1,\dots,N\\
		\partial_tv+\frac{B}{\rho}\partial_x\chi + \frac{1}{\rho}\sum_{\alpha=1}^NA_\alpha^2\partial_x\rho_\alpha +
		v\partial_xv &=0
	\end{align}
\end{subequations}
with the corresponding Jacobian matrix
\begin{equation}\label{jacobi1}
	{\mathbf{A}}(\chi,\rho_1,\dots,\rho_N,v)=
	\begin{pmatrix}
		v & 0 & \dots  & 0 & 0 \\
		0 & v & \dots  & 0 & \rho_1\\
		\vdots &	  & \ddots &   & \vdots\\
		0 & 0 & \dots  & v & \rho_N\\
		\frac{B}{\rho} & \frac{A_1^2}{\rho} && \frac{A_N^2}{\rho} & v
	\end{pmatrix}
\end{equation}

\paragraph{Eigenvalues and eigenvectors.}
We find the eigenvalues of the Jacobian \eqref{jacobi1} that are
\begin{equation}
	\lambda_0=v-A\,,\, \lambda_1=\dots=\lambda_{N}=v\,,\,\lambda_{N+1}=v+A
\end{equation}
where we use the notation
\[
A^2:=\frac{1}{\rho}\sum_{\alpha=1}^NA_\alpha^2\rho_\alpha.
\]
Note that in the single component case $N=1$ we have $A^2=A_1^2$.

In the case $B\ne0$ we obtain the full set of linearly independent eigenvectors
\begin{equation}
	\mathbf{k}_0=\begin{pmatrix}
		0 \\ -\rho_1 \\ \vdots \\ \vdots \\ -\rho_N\\ A
	\end{pmatrix}
	\,,\,
	\mathbf{k}_1=\begin{pmatrix}
		-A_1^2 \\ B \\ 0 \\ \vdots \\ 0 \\ 0
	\end{pmatrix}
	\,,\dots\,,
	\mathbf{k}_N=\begin{pmatrix}
		-A_N^2 \\ 0 \\ \vdots \\ 0 \\ B \\ 0
	\end{pmatrix}
	\,,\,
	\mathbf{k}_{N+1}=\begin{pmatrix}
		0 \\ \rho_1 \\ \vdots \\ \vdots \\ \rho_N\\ A
	\end{pmatrix}
	\,.
\end{equation}

In the case $B=0$ the eigensystem looks different. The reason is that eigenvectors do not depend smoothly on the entry of the matrix. Nevertheless, we have a full system of linearly independent eigenvectors with
\begin{equation}
	\mathbf{k}_1=\begin{pmatrix}
		1 \\ 0 \\ 0\\ \vdots \\ \vdots \\ 0\\ 0
	\end{pmatrix}
	\,,\,
	\mathbf{k}_2=\begin{pmatrix}
		0 \\ -A_2^2 \\ A_1^2 \\ 0\\ \vdots \\ 0 \\ 0
	\end{pmatrix}
	\,,\dots\,,
	\mathbf{k}_N=\begin{pmatrix}
		0 \\ -A_3^2 \\ 0\\ A_1^2 \\ 0 \\ \vdots \\ 0
	\end{pmatrix}
	\,,\,
	\mathbf{k}_{N}=\begin{pmatrix}
		0 \\ -A_N^2 \\ 0\\ \vdots \\ 0 \\ A_1^2 \\ 0
	\end{pmatrix}
\end{equation}
and $\mathbf{k}_0$ and $\mathbf{k}_{N+1}$ as before.

As shown before in both cases we have a full set of linearly independent eigenvector. This implies that system \eqref{hyp} is hyperbolic. For 
the single component case $N=1$ it is even strictly hyperbolic.

\paragraph{Characteristic fields.}

Let $\mathbf{u}=(\chi,\rho_1,\dots,\rho_N,v)^T$ denote the vector of primitive variables.
For $j=0$ or $j=N+1$ we obtain
\[
\nabla_\mathbf{u}\lambda_j \cdot \mathbf{k}_j=A\ne0\,.
\]
This implies that the associated characteristic fields are genuinely non-linear and the corresponding waves are shocks or rarefactions.

For the multiple eigenvalue $\lambda_j=v$ with $j=1,\dots, N$ one can easily verify that 
\[
\nabla_\mathbf{u}\lambda_j \cdot \mathbf{k}_j=0\,.
\]
The associated characteristic field is linearly degenerate and the corresponding wave is a classical contact.

\paragraph{Riemann invariants.}

Let $\lambda$ be an eigenvector of multiplicity $m$ in a system of dimension $n$. Then there exist $n-m$ Riemann invariants 
across the wave corresponding to $\lambda$.

Accordingly we have $N+1$ Riemann invariants across the outer waves belonging to $\lambda_0$ and $\lambda_{N+1}$ and we have $2$ Riemann invariants across the contact wave in the middle.

To find the Riemann invariants across the field $j$, $j=0,\dots,N+1$ one has to solve the system
\[
\frac{\D u_0}{k_{j,0}}=\frac{\D u_1}{k_{j,1}}= \dots = \frac{\D u_{N+1}}{k_{j,{N+1}}}\,.
\]
Consider the \textbf{case $\mathbf{j=0}$}. Then the system of ordinary differential equations to solve becomes
\[
\frac{\D \chi}{0}=\frac{\D \rho_1}{-\rho_1}= \dots = \frac{\D \rho_N}{-\rho_N}=\frac{\D v}{A}\,.
\]
It is easy to see that the phase field $\chi$ is constant across the 0th wave with $\chi\equiv\chi_-$.
For $j=2,\dots,N$ we have
\[
\frac{\D \rho_1}{-\rho_1}=\frac{\D \rho_j}{-\rho_j}
\]
This gives
\[
\ln(\rho_1)-\ln(\rho_j)= const \qquad \mbox{resp.} \qquad \rho_j=c_{j-}\rho_1\,.
\]
It remains to solve
\[
-\frac{A}{\rho_1}\D \rho_1=\D v\,.
\]
Defining $c_{1-}:=1$ we get
\[
\D v=-\frac{1}{\rho_1}\sqrt{\frac{\sum_{\alpha=1}^NA_\alpha^2c_{\alpha-}\rho_1}{\sum_{\alpha=1}^Nc_{\alpha-}\rho_1}} \D \rho_1
= -\frac{1}{\rho_1}\sqrt{\frac{\sum_{\alpha=1}^NA_\alpha^2c_{\alpha-}}{\sum_{\alpha=1}^Nc_{\alpha-}}} \D \rho_1\,.
\]
Keeping in mind that the phase field is an invariant we have $A_\alpha=A_\alpha(\chi_-)=const$ we define $A_{\alpha-}:=A_\alpha(\chi_-)$ and we finally obtain
\[
v=v_- -\ln(\rho_1)\sqrt{\frac{\sum_{\alpha=1}^NA_{\alpha-}^2c_{\alpha-}}{\sum_{\alpha=1}^Nc_{\alpha-}}}
+\ln(\rho_{1-})\sqrt{\frac{\sum_{\alpha=1}^NA_{\alpha-}^2c_{\alpha-}}{\sum_{\alpha=1}^Nc_{\alpha-}}}\,.
\]

The \textbf{case} $\mathbf{j=N+1}$ is quite similar. An analogous calculation with $c_{1+}:=1$ gives the following relations
\begin{align*}
	\chi & \equiv \chi_+\\
	\rho_j & = c_{j+}\rho_1 \quad j=2,\dots,N\\
	v & =v_+ +\ln(\rho_1)\sqrt{\frac{\sum_{\alpha=1}^NA_{\alpha+}^2c_{\alpha+}}{\sum_{\alpha=1}^Nc_{\alpha+}}}
	-\ln(\rho_{1+})\sqrt{\frac{\sum_{\alpha=1}^NA_{\alpha+}^2c_{\alpha+}}{\sum_{\alpha=1}^Nc_{\alpha+}}}
	\qquad \mbox{with $A_{\alpha+}=A_\alpha(\chi_+)$}\,.
\end{align*}

For the \textbf{contact} one immediately can see that the velocity is a Riemann invariant. Nevertheless we fail to determine 
the second invariant. However, from the single component case with pure phases we know, that also the pressure is a constant across the middle wave. To verify if this is true in general we rewrite the quasilinear system, now using the variables $\rho_\alpha$, $\alpha=1,\dots,N$, $v$, $p$.

\paragraph{New choice of variables.}
As mentioned before we rewrite our systems in terms of $\rho_\alpha$, $\alpha=1,\dots,N$, $v$ and $p$. This leads to
\begin{subequations}
	\begin{align}
		\partial_tp+v\partial_xp+\sum_{\alpha=1}^NA^2_\alpha\rho_\alpha\partial_xv &=0\\
		\partial_t\rho_\alpha+v\partial_x\rho_\alpha+\rho_\alpha\partial_xv &=0 \qquad \alpha=1,\dots,N\\
		\partial_tv+\frac{1}{\rho}\partial_xp+v\partial_xv &=0
	\end{align}
\end{subequations}
with the corresponding Jacobian
\[
\begin{pmatrix}
	v & 0 & \dots  & 0 & \sum_{\alpha=1}^NA_\alpha^2\rho_\alpha \\
	0 & v & \dots  & 0 & \rho_1\\
	\vdots &	  & \ddots &   & \vdots\\
	0 & 0 & \dots  & v & \rho_N\\
	\frac{1}{\rho} & 0 & \dots & 0 & v
\end{pmatrix}
\,.
\]
The eigenvectors belonging to the multiple eigenvalue $\lambda_1=\dots=\lambda_N$ are given by
\[
\begin{pmatrix}
	0 \\ 1 \\ 0 \\ \vdots \\ 0 \\ 0
\end{pmatrix}
\,,\,
\dots\,,\,
\begin{pmatrix}
	0 \\ 0 \\ \vdots \\ 0 \\ 1 \\ 0
\end{pmatrix}
\,.
\]
Once again we easily can see that the velocity remains constant across the contact wave. In addition we find the pressure to 
be the further invariant as supposed in the recent paragraph.

%% file: sec4.tex
\section{Exact solution of the Riemann problem}
In this section we want to construct the exact solution to the Riemann problem \eqref{hyp}, \eqref{pressure}, \eqref{ind}. For that we summarize the findings of the last section.
\begin{itemize}
	\item
	The solution of the Riemann problem consists of four constant states that are separated by three waves. The middle wave is a contact while the outer waves are shocks or rarefactions.
	\item 
	The phase field $\chi$ may change across the contact wave, but stays constant everywhere else. This means that the initial profile of the phase field is shifted with the flow.
	\item
	Due to the fact that the solution for the phase field is known, it remains to solve the system consisting of the partial mass balances and the momentum balance. This system is in divergence form. Accordingly the Rankine-Hugeniot jump conditions are satisfied across discontinuities. These are given by
	\begin{align}
		\label{rhmb}
		\rho_\alpha''v''-\rho_\alpha'v' &= S(\rho_\alpha''-\rho_\alpha') \qquad \alpha=1,\dots,N \\
		\label{rhib}
		\sum_{\alpha=1}^N (\rho_\alpha''(v'')^2-\rho_\alpha'(v')^2)+(p''-p') &= S\sum_{\alpha=1}^N(\rho_\alpha''v''-\rho_\alpha'v')
	\end{align}
	where $S$ denotes the propagation speed of the discontinuity and $'$ and $''$ indicate the states to the left and to the right of the discontinuity, resp.
	\item
	The velocity and the pressure are Riemann invariants across the contact wave. This allows to follow the strategy described in the book of Toro \cite{toro} to construct the Riemann solution.
\end{itemize}

\begin{figure}[h!]
	\begin{center}
		\begin{tikzpicture}[scale=1.5]
			\draw [dotted, thick] (0,0) -- (0.5,3);
			\draw [thick] (0,0) -- (-3,3);
			\draw [thick] (0,0) -- (-3.2,3);
			\draw [thick] (0,0) -- (3,3);
			\draw [thick] (0,0) -- (3.2,3);
			\draw [>=latex,->] (-3.7,0) -- (3,0) node[pos=.55, below]{0} node[pos=1, below]{$x$};
			\draw [>=latex,->] (-3.7,0) -- (-3.7,4) node[pos=1, left]{$t$};
			\node at (-2.5,1) {{\textit{$\chi_-,\rho_{\alpha-},v_-,p_-$}}};
			\node at (-1.0,2.3) {{\textit{$\chi_-,\rho_{\alpha-}^*,v^*,p^*$}}};
			\node at (1.25,2.3) {{\textit{$\chi_+,\rho_{\alpha+}^*,v^*,p^*$}}};
			\node at (2.3,1) {{\textit{$\chi_+,\rho_{\alpha+},v_+,p_+$}}};
		\end{tikzpicture}
		
	\end{center}
	\caption{Structure of the Riemann solution}
	\label{wavepattern}
\end{figure}
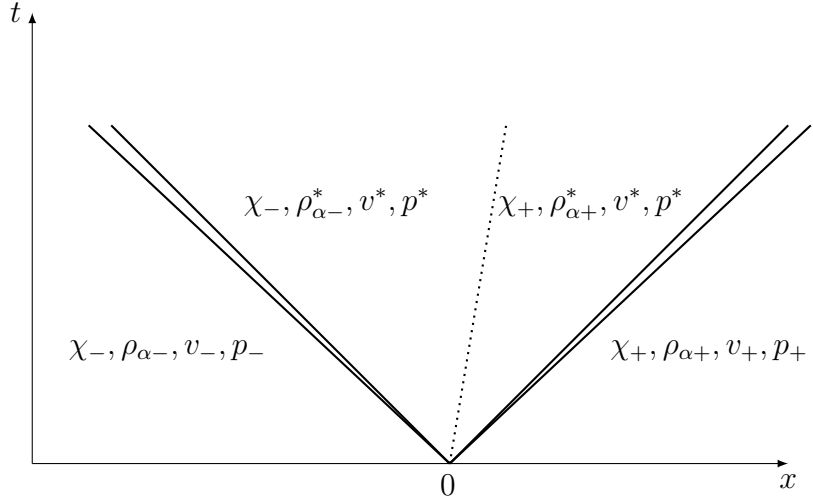

\paragraph{Solution strategy.}
The structure of the Riemann solution is depicted in Figure \ref{wavepattern}.
The region between the two outer wave is called \textit{star region}. This region is divided into two subregions \textit{star left} and \textit{star right} with the corresponding vectors of unknowns $\mathbf{W}_-^*$ and $\mathbf{W}_+^*$. Using $\chi_\pm^*=\chi_\pm$ and $v_-^*=v_+^*=:v^*$ these vectors read
$\mathbf{W}_\pm^*=(\chi_\pm,\rho_{1\pm}^*,\dots,\rho_{N\pm}^*,v^*)^T$.

To determine $\mathbf{W}_\pm^*$ we take advantage of $v_-^*=v_+^*=:v^*$ and $p_-^*=p_+^*=:p^*$.
Our aim is to derive a function
\begin{align}
	f(p,\mathbf{W}_-,\mathbf{W}_+)=f_-(p,\mathbf{W_-})+f_+(p,\mathbf{W_+})+(v_+-v_-)
\end{align}
such that the only root $p=p_*$ is the solution for the pressure in the star region.
Here the functions $f_-$ and $f_+$ relate the initial velocities $v_-$ and $v_+$ to $v^*$ only in terms of the initial data $\mathbf{W}_\pm$ and the unknown pressure $p^*$, i.e.
\begin{align}
	v^*=v_--f_-(p^*,\mathbf{W}_-) \qquad \mbox{and} \qquad v^*=v_++f_+(p^*,\mathbf{W}_+)\,.
\end{align}
This construction uses the fact, that the pressure is a constant in the star region. Accordingly we choose $p^*$ to be the unknown and we have to eliminate all partial densities $\rho_{\alpha\pm}$.

\paragraph{Rarefactions.}

Assume the left wave is a rarefaction. Then from the last section we know that
\begin{equation}\label{vst1}
	v^*=v_- -\ln(\rho_{1-}^*)\sqrt{\frac{\sum_{\alpha=1}^NA_{\alpha-}^2c_{\alpha-}}{\sum_{\alpha=1}^Nc_{\alpha-}}}
	+\ln(\rho_{1-})\sqrt{\frac{\sum_{\alpha=1}^NA_{\alpha-}^2c_{\alpha-}}{\sum_{\alpha=1}^Nc_{\alpha-}}}
\end{equation}
and $\rho_{j-}^*=c_{j-}\rho_{1-}^*$ with $c_{j-}=\frac{\rho_{j-}}{\rho_{1-}}$ and $j=1,\dots,N$.
Using \eqref{pressure} we find
\begin{align*}
	p^* &= -W(\chi_-)+h(\chi_-)\sum_{\alpha=1}^N(a_{L\alpha}^2 c_{j-}\rho_{1-}^*+d_{L\alpha}) +
	(1-h(\chi_-))\sum_{\alpha=1}^Na_{V\alpha}^2c_{j-}\rho_{1-}^*\\
	&= {\cal{A}}_{0-}+{\cal{A}}_{1-}\rho_{1-}^*
\end{align*}
with
\begin{equation}\label{defA}
	{\cal{A}}_{0-}=-W(\chi_-)+h(\chi_-)\sum_{\alpha=1}^Nd_{L\alpha}
	\qquad \mbox{and} \qquad
	{\cal{A}}_{1-}=h(\chi_-)\sum_{\alpha=1}^Na_{L\alpha}^2 c_{j-}
	+(1-h(\chi_-))\sum_{\alpha=1}^Na_{V\alpha}^2c_{j-}\,.
\end{equation}
So we can replace $\rho_{1-}^*$ in \eqref{vst1} by $\rho_{1-}^*=\frac{p^*-{\cal{A}}_{0-}}{{\cal{A}}_{1-}}$.

An analogous calculation for a right rarefaction leads to
\begin{equation}
	v^*  = v_+ +\ln\left(\frac{p^*-{\cal{A}}_{0+}}{{\cal{A}}_{1+}}\right)\sqrt{\frac{\sum_{\alpha=1}^NA_{\alpha+}^2c_{\alpha+}}{\sum_{\alpha=1}^Nc_{\alpha+}}}
	-\ln(\rho_{1+})\sqrt{\frac{\sum_{\alpha=1}^NA_{\alpha+}^2c_{\alpha+}}{\sum_{\alpha=1}^Nc_{\alpha+}}}
	\qquad \mbox{with $A_{\alpha+}=A_\alpha(\chi_+)$}\,.
\end{equation}

\paragraph{Shocks.}
Now let us assume the left wave is a shock propagating with speed $S_-$. We define the left relative mass fluxes $Q_{\alpha-}$ be rewriting \eqref{rhmb} as
\begin{align}
	\label{relmb}
	\rho_{\alpha-}(v_--S_-)=\rho_{\alpha-}^*(v^*-S_-)\quad=:Q_{\alpha-} \qquad \alpha=1,\dots,N\,.
\end{align}
We note that 
\begin{align}
	\frac{\rho_{\alpha-}^*}{\rho_{\alpha-}}=\frac{v_--S_-}{v^*-S_-} \qquad \mbox{for all $\alpha=1,\dots,N$.}
\end{align}
Accordingly 
\begin{align}
	\rho_{\alpha_-}^*=\rho_{\alpha-}\frac{\rho_{1-}^*}{\rho_{1-}} \qquad \mbox{for all $\alpha=1,\dots,N$.}
\end{align}
Using \eqref{relmb} this we can rewrite \eqref{rhib} as
\begin{align}\label{relib}
	(v_--v^*)\sum_{\alpha=1}^NQ_{\alpha-}+(p_--p^*)=0
\end{align}
and we obtain 
\begin{align}
	v^*=v_--\frac{p^*-p_-}{\sum_{\alpha=1}^NQ_{\alpha-}}
\end{align}
where we will express $\sum_{\alpha=1}^NQ_{\alpha-}$ in terms of $p^*$. To do this we rewrite \eqref{relib} as follows.
\begin{align}\label{sumQ}
	\sum_{\alpha=1}^NQ_{\alpha-}=-\frac{p_--p^*}{v_--v^*}=-\frac{p_--p^*}{(v_--S_-)-(v^*-S_-)}
	=-\frac{p_--p^*}{\frac{Q_{j-}}{\rho_{j-}}-\frac{Q_{j-}}{\rho_{j-}^*}}\,.
\end{align}
Multiplying \eqref{sumQ} by $Q_{j-}$ and summing up over $j$ we get
\begin{align}
	\sum_{j=1}^NQ_{j-}\sum_{\alpha=1}^NQ_{\alpha-}=-\sum_{j=1}^N
	\frac{p_--p^*}{\frac{1}{\rho_{j-}}-\frac{1}{\rho_{j-}^*}}\,.
\end{align}
This leads to
\begin{align}
	\sum_{\alpha=1}^NQ_{\alpha-}=\sqrt{-\sum_{\alpha=1}^N\frac{p_--p^*}{\frac{1}{\rho_{\alpha-}}-\frac{1}{\rho_{\alpha-}^*}}}
\end{align}
and we finally obtain
\begin{align}
	v^*=v_--\frac{\sqrt{p^*-p_-}}{\sqrt{\sum_{\alpha=1}^N\frac{\rho_{\alpha-}\rho_{\alpha-}^*}{\rho_{\alpha-}^*-\rho_{\alpha-}}}}
	=v_--\frac{\sqrt{p^*-p_-}}{\sqrt{\frac{\rho_{1-}^*}{\rho_{1-}^*-\rho_{1-}}\sum_{\alpha=1}^N\rho_{\alpha-}}}
	=v_--\frac{\sqrt{p^*-p_-}}{\sum_{\alpha=1}^N\rho_{\alpha-}}\sqrt{1-\frac{{\cal{A}}_{1-}\rho_{1-}}{p^*-{\cal{A}}_{0-}}}
	\,.
\end{align}
Analogously we obtain for a right shock
\begin{align}
	v^*=v_++
	\frac{\sqrt{p^*-p_+}}{\sum_{\alpha=1}^N\rho_{\alpha+}}\sqrt{1-\frac{{\cal{A}}_{1+}\rho_{1+}}{p^*-{\cal{A}}_{0+}}}\,.
\end{align}

\paragraph{Existence and uniqueness of the Riemann solution.}
We summarize the results of the last paragraphs in the following
\begin{teo}\label{teo1}
	Let be given the function
	\[
	f(p,\mathbf{W}_-,\mathbf{W}_+)=f_-(p,\mathbf{W_-})+f_+(p,\mathbf{W_+})+(v_+-v_-)
	\]
	with
	\begin{align*}
		f_-(p,\mathbf{W_-})=
		\begin{cases}
			\frac{\sqrt{p-p_-}}{\sum_{\alpha=1}^N\rho_{\alpha-}}\sqrt{1-\frac{{\cal{A}}_{1-}\rho_{1-}}{p-{\cal{A}}_{0-}}}
			& \mbox{if $p>p_-$ (shock)}\\
			\ln\left(\frac{p-{\cal{A}}_{0-}}{{\cal{A}}_{1-}}\right)\sqrt{\frac{\sum_{\alpha=1}^NA_{\alpha-}^2c_{\alpha-}}{\sum_{\alpha=1}^Nc_{\alpha-}}}
			-\ln(\rho_{1-})\sqrt{\frac{\sum_{\alpha=1}^NA_{\alpha-}^2c_{\alpha-}}{\sum_{\alpha=1}^Nc_{\alpha-}}}
			& \mbox{if $p\le p_-$ (rarefaction)}
		\end{cases}
		\\
		f_+(p,\mathbf{W_+})=
		\begin{cases}
			\frac{\sqrt{p-p_+}}{\sum_{\alpha=1}^N\rho_{\alpha+}}\sqrt{1-\frac{{\cal{A}}_{1+}\rho_{1+}}{p-{\cal{A}}_{0+}}}
			& \mbox{if $p>p_+$ (shock)}\\
			\ln\left(\frac{p-{\cal{A}}_{0+}}{{\cal{A}}_{1+}}\right)\sqrt{\frac{\sum_{\alpha=1}^NA_{\alpha+}^2c_{\alpha+}}{\sum_{\alpha=1}^Nc_{\alpha+}}}
			-\ln(\rho_{1+})\sqrt{\frac{\sum_{\alpha=1}^NA_{\alpha+}^2c_{\alpha+}}{\sum_{\alpha=1}^Nc_{\alpha+}}}
			& \mbox{if $p\le p_+$ (rarefaction)}
		\end{cases}
	\end{align*}	
	with
	\begin{align*}
		{\cal{A}}_{0\pm} &= -W(\chi_\pm)+h(\chi_\pm)\sum_{\alpha=1}^Nd_{L\alpha}
		\\
		{\cal{A}}_{1\pm} &= h(\chi_\pm)\sum_{\alpha=1}^Na_{L\alpha}^2 c_{j\pm}
		+(1-h(\chi_\pm))\sum_{\alpha=1}^Na_{V\alpha}^2c_{j\pm}
		\\
		c_{\alpha\pm} &= \frac{\rho_{\alpha\pm}}{\rho_{1\pm}} \qquad \alpha=1,\dots,N.
	\end{align*}
	Then the function $f(p,\mathbf{W}_-,\mathbf{W}_+)$ has a unique root $p=p^*$
	that is the unique solution for the pressure $p^*$ of the Riemann problem \eqref{hyp}, \eqref{pressure}, \eqref{ind}. The velocity $v^*$ can be calculated using
	\[
	v^*=\frac{1}{2}(v_-+v_+)+\frac{1}{2}(f_+(p^*,\mathbf{W_+}) - f_-(p^*,\mathbf{W_-}))\,.
	\]	
\end{teo}

\begin{proof}
	The function $f$ is strictly increasing in $p$. For $p\to\min\{{\cal{A}}_{0-},{\cal{A}}_{0+}\}$ the function $f$ tends to $-\infty$. For $p\to+\infty$ we have $f\to+\infty$. Accordingly $f$ has a unique root that by construction is the solution for the pressure $p^*$ of the problem considered. The remaining part of the theorem is obvious.
\end{proof}

\begin{rem}
	To determine the remaining unknown quantities of the solution one has to use the relations above.
	Here one has to take care for the type of the waves.
\end{rem}

\begin{rem}
	For the special case $N=1$, $\chi_-=-1$ and $\chi_+=1$ Theorem \ref{teo1} reduces to Theorem 6.2 (Solution of isothermal two-phase Euler equations without phase transition) in \cite{hdw}.
\end{rem}

%% file: sec5.tex
\section{Examples and Benchmarks}

In the following we will give the initial data and the full solution for two examples. These examples may be used as benchmarks to validate some numerical methods.

\paragraph{Example 1.}
First we consider a 2-component example. The initial data and parameters used are summarized in Table \ref{c1in}.
\begin{table}[h!]
	\begin{center}
		\begin{tabular}{|c||c|c|c|c||c||c|c|c|c|}
			\hline& $\chi$ & $\rho_1$ & $\rho_2$ & $v$ &  & $a_1$ & $a_2$ & $d_1$ & $d_2$ \\
			\hline
			Left & $-0.95$ & $2.5$ & $7.5$ & $0$ & Vapor & $200$ & $300$ & $0$ & $0$\\
			\hline
			Right & $0.5$ & $ 600$ & $800$ & $0$ & Liquid & $500$ & $400$ & $-1.495\cdot10^8$ & $-6.35\cdot10^7$\\
			\hline
		\end{tabular}
	\end{center}
	\caption{Initial data and parameters Example 1}
	\label{c1in}
\end{table}	

\begin{figure}[h!]
	\includegraphics{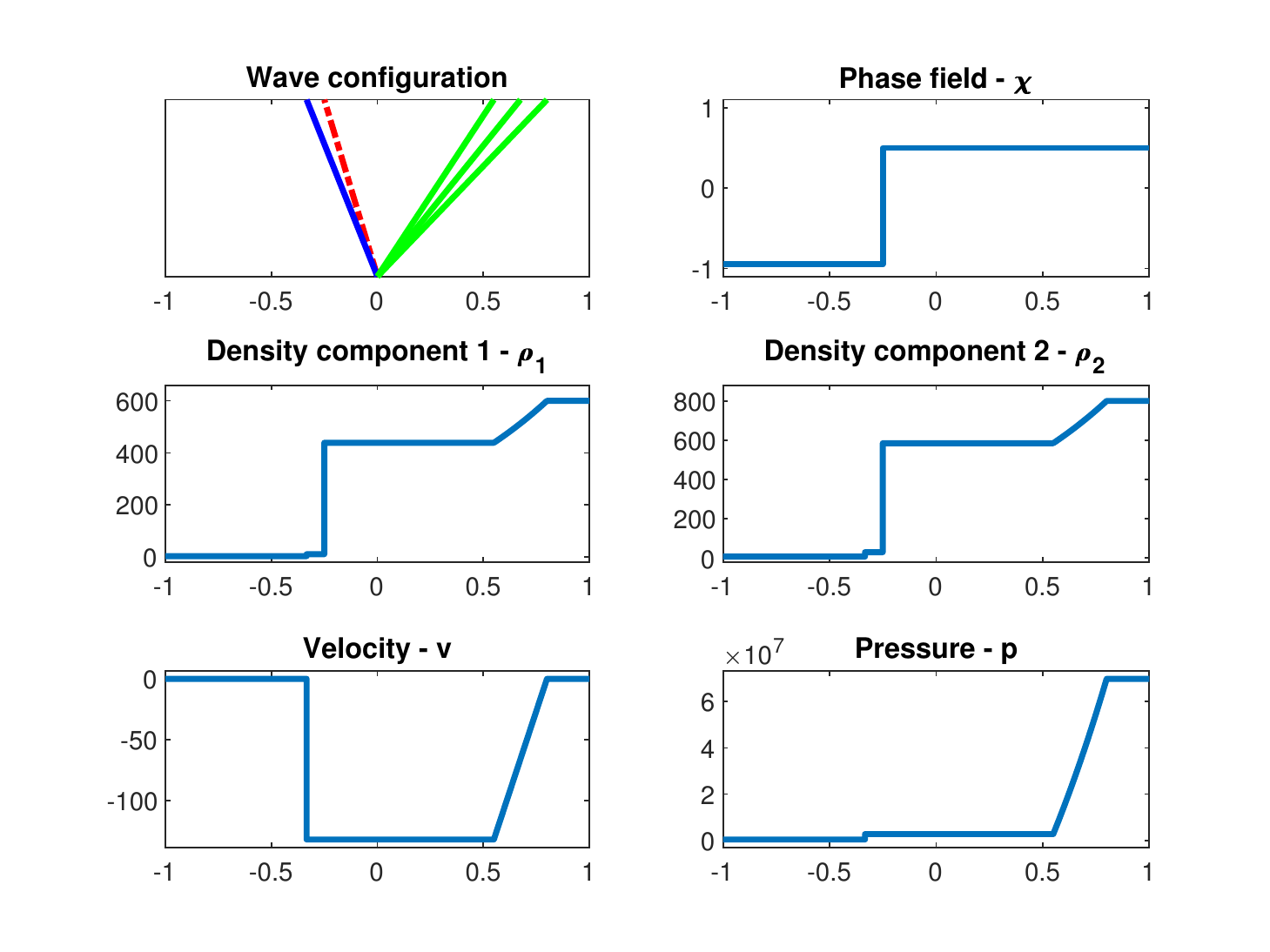}
	\caption{Exact solution Example 1}
	\label{fex1}
\end{figure}

The solution consists of 4 constant states, separated by a left shock, a contact discontinuity and a right rarefaction, see Figure \ref{fex1}.

For the wave speeds and the states in the star region we obtain

\begin{table}[h!]
	\begin{center}
		\begin{tabular}{|c|c|c|c|c|}
			\hline 
			$S_{L}$ & $S_{R,tail}$ & $S_{R,head}$ & $p^*$ & $v^*$ \\
			\hline
			$-176.3412$ & $289.925$ & $422.207$ & $2716903.0964$ & $-132.2825$\\
			\hline
		\end{tabular}
	\end{center}
	\caption{Solution Example 1}
	\label{c1loes}
\end{table}	

\paragraph{Example 2.}
Next we consider a 3-component example. The initial data and parameters used are summarized in Tables \ref{c2ina} and \ref{c2inb}.
\begin{table}[h!]
	\begin{center}
		\begin{tabular}{|c||c|c|c|c|c|}
			\hline& $\chi$ & $\rho_1$ & $\rho_2$ & $\rho_3$ & $v$ \\
			\hline
			Left & $-0.95$ & $2.5$ & $7.5$ & $1$ & $-50$ \\
			\hline
			Right & $0.5$ & $ 300$ & $800$ & $250$ & $20$ \\
			\hline
		\end{tabular}
	\end{center}
	\caption{Initial data Example 2}
	\label{c2ina}
\end{table}	

\begin{table}[h!]
	\begin{center}
		\begin{tabular}{|c||c|c|c|c|c|c|}
			\hline  & $a_1$ & $a_2$ & $a_3$ & $d_1$ & $d_2$ & $d_3$\\
			\hline
			 Vapor & $200$ & $300$ & $100$ & $0$ & $0$ & $0$\\
			\hline
			 Liquid & $250$ & $400$ & $200$ & $-7.45\cdot10^7$ & $-6.35\cdot10^7$ & $-3.15\cdot10^7$ \\
			\hline
		\end{tabular}
	\end{center}
	\caption{Parameters Example 2}
	\label{c2inb}
\end{table}	

\begin{figure}[h!]
	\includegraphics{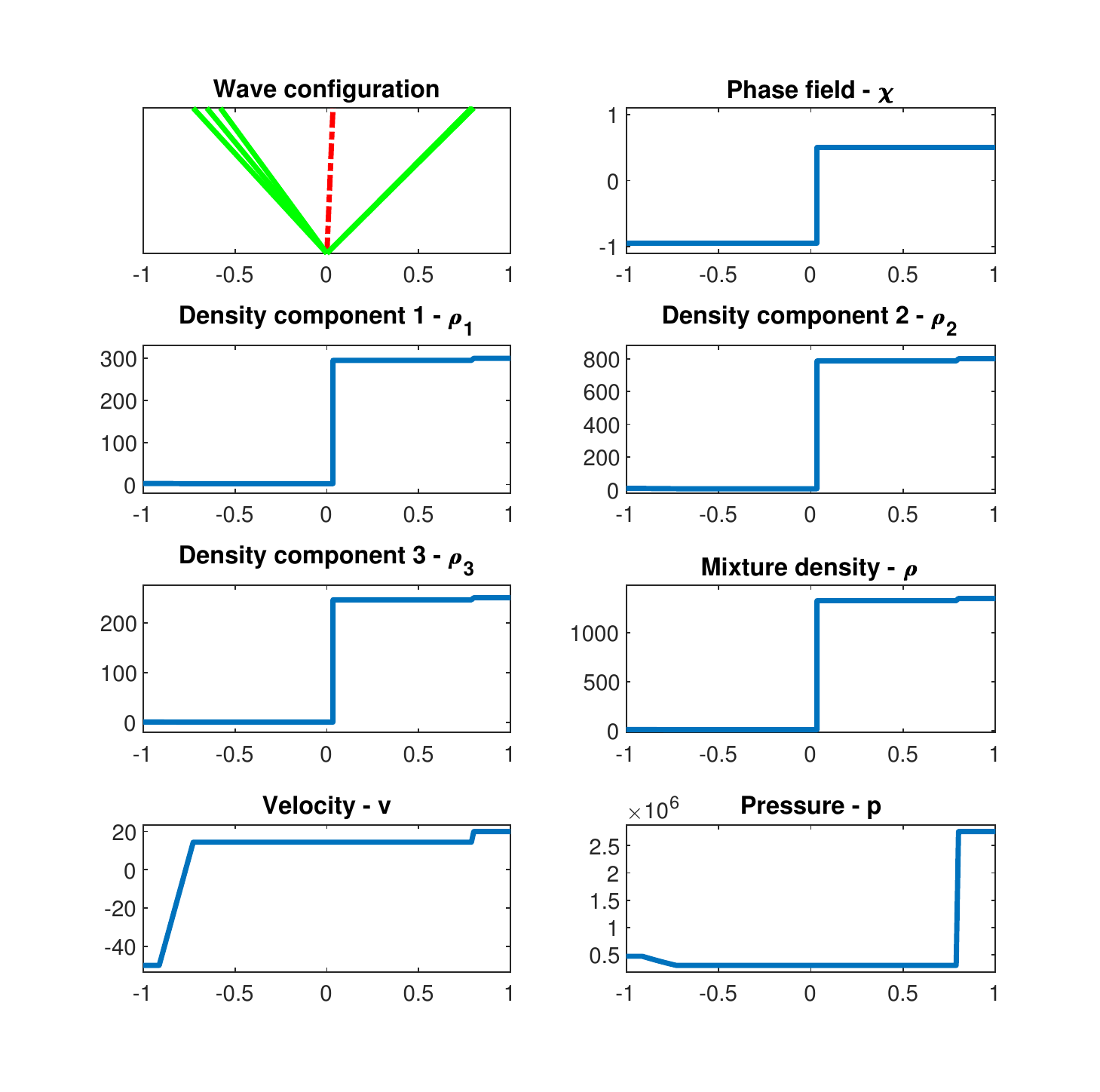}
	\caption{Exact solution Example 2}
	\label{fex2}
\end{figure}

The solution consists of 4 constant states, separated by a left rarefaction, a contact discontinuity and a right rarefaction, see Figure \ref{fex2}.

For the wave speeds and the states in the star region we obtain

\begin{table}[h!]
	\begin{center}
		\begin{tabular}{|c|c|c|c|c|c|}
			\hline 
			$S_{L,head}$ & $S_{L,tail}$ & $S_{R,tail}$ & $S_{R,head}$ & $p^*$ & $v^*$ \\
			\hline
			$-317.331$ & $-252.907$ & $343.028$ & $348.604$ & $305261.3806$ & $14.4244$\\
			\hline
		\end{tabular}
	\end{center}
	\caption{Solution Example 2}
	\label{c2loes}
\end{table}	

%% file: hantke_mwy_ana.bbl
\begin{thebibliography}{10}

\bibitem{baernun}
M.~Baer and J.~Nunziato.
\newblock A two-phase mixture theory for the deflagration-to-detonation
  transition ({D}{D}{T}) in reactive granular materials.
\newblock {\em International Journal of Multiphase Flow}, 12(6):861 -- 889,
  1986.

\bibitem{chalons1}
C.~Chalons, P.~Engel, and C.~Rohde.
\newblock A conservative and convergent scheme for undercompressive shock
  waves.
\newblock {\em SIAM Journal on Numerical Analysis}, 52(1):554--579, 2014.

\bibitem{chalons2}
C.~Chalons, C.~Rohde, and M.~Wiebe.
\newblock A finite volume method for undercompressive shock waves in two space
  dimensions.
\newblock {\em {ESAIM: Mathematical Modelling and Numerical Analysis}},
  51(5):1987--2015, 2017.

\bibitem{db}
W.~Dreyer and D.~{Bothe}.
\newblock Continuum thermodynamics of chemically reacting fluid mixtures.
\newblock {\em Acta Mechanica}, 226:1757--1805, 2015.

\bibitem{dgk}
W.~Dreyer, J.~Giesselmann, and C.~{Kraus}.
\newblock A compressible mixture model with phase transition.
\newblock {\em Physics D: Nonlinear Phenomena}, 273:1--13, 2014.

\bibitem{fechter}
S.~Fechter, C.-D. Munz, C.~Rohde, and C.~{Zeiler}.
\newblock A sharp interface method for compressible liquid -- vapor flow with
  phase transition and surface tension.
\newblock {\em Journal of Computational Physics}, 336:347--374, 2017.

\bibitem{hdw}
M.~Hantke, W.~Dreyer, and G.~{Warnecke}.
\newblock Exact solutions to the {R}iemann problem for compressible isothermal
  {E}uler equations for two phase flows with and without phase transitions.
\newblock {\em Quarterly of Applied Mathematics}, LXXI 3:509--540, 2013.

\bibitem{hmwy}
M.~Hantke, C.~Matern, G.~Warnecke, and H.~Yaghi.
\newblock A new method to discretize a multi component phase field model.
\newblock To appear.

\bibitem{hm}
M.~Hantke and S.~M\"uller.
\newblock {A}nalysis and simulation of a new multi-component two-phase flow
  model with phase transitions and chemical reactions.
\newblock {\em Quarterly of applied mathematics}, 76(2):253--287, 2018.

\bibitem{herard}
J.~H\`{e}rard.
\newblock A three phase flow model.
\newblock {\em Math. Comput. Model}, 45:732--755, 2007.

\bibitem{simapa}
S.~M{\"u}ller, M.~Hantke, and P.~{Richter}.
\newblock Closure conditions for non-equilibrium multi-component models.
\newblock {\em Continuum mechanics and thermodynamics}, 28:1157--1189, 2016.

\bibitem{romenski}
E.~Romenski, A.~D. Resnyansky, and E.~F. Toro.
\newblock Conservative hyperbolic formulation for compressible two-phase flow
  with different phase pressures and temperatures.
\newblock {\em Quarterly of Applied Mathematics}, 65:259--79, 2007.

\bibitem{schleper}
V.~Schleper.
\newblock A hll-type riemann solver for two-phase flow with surface forces and
  phase transitions.
\newblock {\em Applied Numerical Mathematics}, 108:256--270, 2016.

\bibitem{theinRD}
F.~Thein, E.~Romenski, and M.~Dumbser.
\newblock Exact and numerical solutions of the riemann problem for a
  conservative model of compressible two-phase flows, 2022.

\bibitem{toro}
E.~Toro.
\newblock {\em Riemann Solvers and Numerical Methods for Fluid Dynamics}.
\newblock Springer-Verlag, Berlin, third edition, 2009.
\newblock A practical introduction.

\end{thebibliography}
